# Solving the Fermat and Fibonacci Equations with the Lambert-Tsallis Wq Function

R. V. Ramos

rubens.ramos@ufc.br

*Lab. of Quantum Information Technology, Department of Teleinformatic Engineering – Federal University of Ceara - DETI/UFC, C.P. 6007 – Campus do Pici - 60455-970 Fortaleza-Ce, Brazil.*

*Abstract* — In this work, the Lambert-Tsallis $W_q$ function is used to provide analytical solutions of fractional polynomials of the type $ax^r+bx^s+c = 0$. This class of fractional polynomial appears in several areas of physics as well it is in the heart of some famous mathematical problems, like the Fermat and Fibonacci's equations. Therefore, analytical solutions for the equations $A^x + B^x = C^x$ and $\varphi^x - \bar\varphi^x = y\sqrt{5}$, where $\varphi = (1+\sqrt{5})/2$ and $\bar\varphi = (\sqrt{5}-1)/2$, are also provided.

*Keywords* — **Lambert-Tsallis $W_q$ function, Fractional polynomial, Fermat equation, Fibonacci equation.**

## 1. Introduction

There is a large number of famous equations whose analytical solutions are not known. Their solutions could give an impulse in the development of different areas of mathematics and physics. Basically, this happens because the required function able to solve those equations was not yet developed and, in the best of my knowledge, they cannot be solved by the set of elementary functions currently available. Two examples of such equations are: I) the Fermat equation: $A^x + B^x = C^x$ and II) the Fibonacci equation $\varphi^x - \bar\varphi^x = y\sqrt{5}$, where $\varphi = (1+\sqrt{5})/2$ and $\bar\varphi = (\sqrt{5}-1)/2$. It is a trivial mathematical exercise to show that these two equations are particular cases of the following fractional polynomial equation

$$ax^\alpha + bx^\beta + c = 0. \tag{1}$$

In this direction, by using the Lambert-Tsallis $W_q$ function [1], this work provides the analytical solutions of eq. (1) that is, $x(a, \alpha, b, \beta, c)$, and, by consequence, the analytical solutions of Fermat and Fibonacci's equations.

The present work is outlined as follows: In Section 2 the Lambert-Tsallis $W_q$ function is briefly reviewed. In Section 3 the solutions of eq. (1) as well the solutions of the particular cases I) and II) are provided. At last, the conclusions are drawn in Section IV.

## 2. The Lambert-Tsallis Wq function

The Lambert *W* function is an important elementary mathematical function that finds applications in different areas of mathematics and physics [2-5]. The Lambert *W* function is defined as the solution of the equation

$$W(z)e^{W(z)} = z. \tag{2}$$

It is well known that eq. (2) has infinite solutions, nonetheless, only two of them return a real value when the argument *z* is real. In the interval $-1/e \leq z \leq 0$ there exist two real values of $W(z)$: the branch for which $W(z) \geq -1$ is the principal branch named $W_0(z)$ while the branch satisfying $W(z) \leq -1$ is named $W_{-1}(z)$. For $z \geq 0$ only $W_0(z)$ is real and for $z < -1/e$ real solutions do not exist. The point $(z_b = -1/e, W(z_b) = -1)$ is the branch point where the solutions $W_0$ and $W_{-1}$ have the same value and $dW/dz = \infty$. The plot of $W(z)$ versus *z* is shown in Fig. 1.

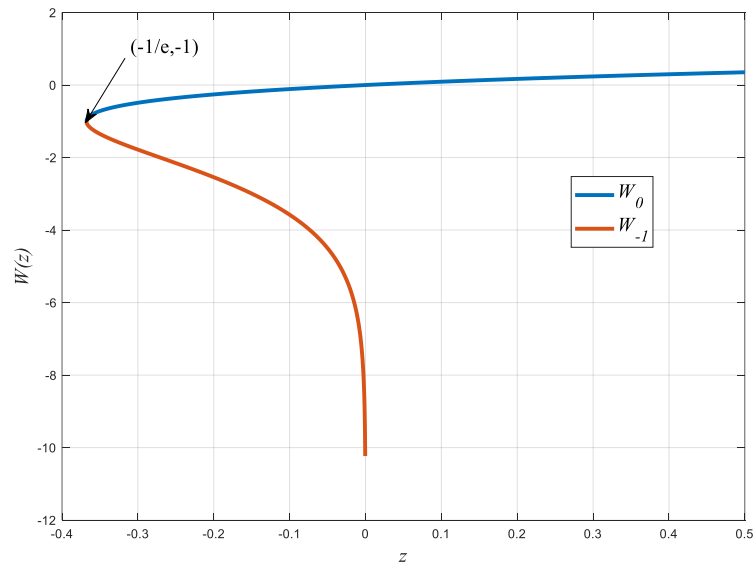

Fig. 1 – *W*(*z*) versus *z*

There are different generalizations of the Lambert $W$ function. A very useful one is the generalization of $W(z)$ that can be obtained when instead of the exponential function in eq. (2) one uses the generalized Tsallis $q$-exponential: $\exp_q(z) = [1+(1-q)z]^{1/(1-q)}$ for $q \neq 1$ [6]. Thus, one gets the Lambert-Tsallis equation,

$$W_q(z)e_q^{W_q(z)} = z, \tag{3}$$

whose solution is the Lambert-Tsallis $W_q$ function introduced in [1]. For $q = 1$ one has $\exp_{q=1}(z) = e^z$ and $W_1(z) = W(z)$. The branch point of $W_q(z)$ is $(z_b = \exp_q(1/(q-2))/(q-2), W_q(z_b) = 1/(q-2))$. One can find the analytical expressions for $W_q(z)$ when $q$ assumes some special values [1]. For example, for $q$ equal to 1/2, 3/2, 4/3 and 2, one has the following upper branches of the respective Lambert-Tsallis functions ($W_{q=a}(z) = W_a(z)$):

$$W_{1/2}(z) = \frac{\left[3\sqrt[3]{2z+\sqrt{\left(2z+\frac{8}{27}\right)^2 - \frac{64}{729}} + \frac{8}{27}} - 2\right]^2}{9\sqrt[3]{2z+\sqrt{\left(2z+\frac{8}{27}\right)^2 - \frac{64}{729}} + \frac{8}{27}}}; \quad z \in (-0.29629, \infty) \tag{4}$$

$$W_{3/2}(z) = \frac{2(z+1) - 2\sqrt{2z+1}}{z}; \quad z \in (-1/2, \infty) \tag{5}$$

$$W_{4/3}(z) = \sqrt[3]{\frac{1}{2}\sqrt[3]{\frac{6561z+2916}{z^3} - \frac{81}{z}} - \frac{9\sqrt[3]{2}}{z\sqrt[3]{\sqrt{\frac{6561z+2916}{z^3} - \frac{81}{z}}}} + 3}; \quad z \in (0, \infty) \tag{6}$$

$$W_2(z) = \frac{z}{z+1}; \quad z \in (-1, \infty). \tag{7}$$

As one may note, $W_2(z)$ does not have a finite branch point. For a general value of $q$ one can calculate the value of $W_q(z)$ numerically (one may note that functions like $\log(z)$, $\exp(z)$, $\tan(z)$, $W(z)$, among others are also calculated numerically). Figure 2 shows the plots of $W_{3/4}(z)$ and $W_{5/3}(z)$ versus $z$.

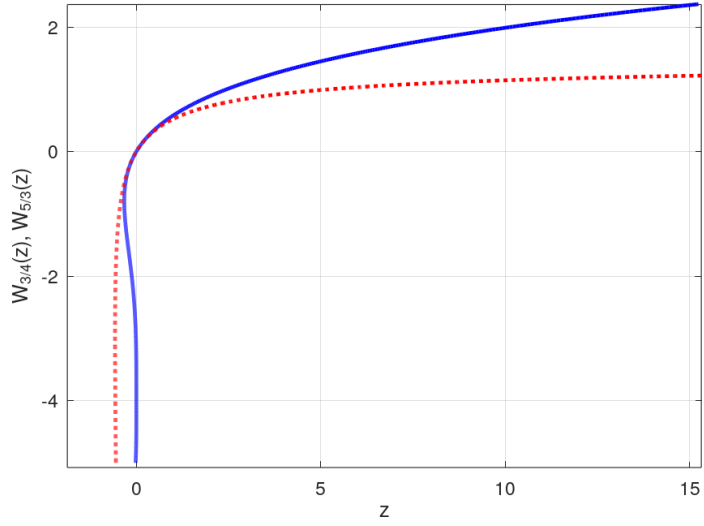

Fig. 2 – $W_q(z)$ versus $z$. $q = 3/4$ (continuous line), $q = 5/3$ (dotted line).

More details about the properties of the Lambert-Tsallis $W_q$ function, its numerical calculation and applications in different areas of physics and engineering can be found in the references [7-14].

3. **The solutions of $ax^\alpha + bx^\beta + c = 0$**

The solutions of $ax^\alpha + bx^\beta + c = 0$ are

$$x_1 = \left[\frac{a}{b}\left(\frac{\alpha}{\beta-\alpha}\right)W_{1-\frac{\alpha}{\beta-\alpha}}\left(\frac{b}{a}\left(\frac{\beta-\alpha}{\alpha}\right)\left(-\frac{c}{a}\right)^{\frac{\beta-\alpha}{\alpha}}\right)\right]^{\frac{1}{\beta-\alpha}} \qquad (8)$$

$$x_2 = \left[\frac{b}{a}\left(\frac{\beta}{\alpha-\beta}\right)W_{1-\frac{\beta}{\alpha-\beta}}\left(\frac{a}{b}\left(\frac{\alpha-\beta}{\beta}\right)\left(-\frac{c}{b}\right)^{\frac{\alpha-\beta}{\beta}}\right)\right]^{\frac{1}{\alpha-\beta}} \qquad (9)$$

The $W_q$ function is multivalued function, hence eqs. (8) and (9) can represent a set of solutions when the argument of $W_q$ is complex. In order to get the solutions of eq. (1) using the Lambert-Tsallis $W_q$ function one has to write eq. (1) in the form shown in eq. (3). This can be done in the following way

$$ax^\alpha + bx^\beta = -c \qquad (10.1)$$

$$ax^\alpha\left(1+\frac{b}{a}x^{\beta-\alpha}\right)=-c \Rightarrow x^\alpha e_0^{\frac{b}{a}x^{\beta-\alpha}} = -\frac{c}{a} \qquad (10.2)$$

$$\left(x^\alpha e_0^{\frac{b}{a}x^{\beta-\alpha}}\right)^{\frac{\beta-\alpha}{\alpha}} = \left(-\frac{c}{a}\right)^{\frac{\beta-\alpha}{\alpha}} \Rightarrow x^{\beta-\alpha} e_{1-\frac{\alpha}{\beta-\alpha}}^{\frac{b}{a}\frac{\beta-\alpha}{\alpha}x^{\beta-\alpha}} = \left(-\frac{c}{a}\right)^{\frac{\beta-\alpha}{\alpha}} \qquad (10.3)$$

$$\frac{b}{a}\frac{\beta-\alpha}{\alpha}x^{\beta-\alpha}e_{1-\frac{\alpha}{\beta-\alpha}}^{\frac{b}{a}\frac{\beta-\alpha}{\alpha}x^{\beta-\alpha}} = \frac{b}{a}\frac{\beta-\alpha}{\alpha}\left(-\frac{c}{a}\right)^{\frac{\beta-\alpha}{\alpha}} \Rightarrow \frac{b}{a}\frac{\beta-\alpha}{\alpha}x^{\beta-\alpha} = W_{1-\frac{\alpha}{\beta-\alpha}}\left[\frac{b}{a}\frac{\beta-\alpha}{\alpha}\left(-\frac{c}{a}\right)^{\frac{\beta-\alpha}{\alpha}}\right] \qquad (10.4)$$

$$x = \left\{\frac{a}{b}\frac{\alpha}{\beta-\alpha}W_{1-\frac{\alpha}{\beta-\alpha}}\left[\frac{b}{a}\frac{\beta-\alpha}{\alpha}\left(-\frac{c}{a}\right)^{\frac{\beta-\alpha}{\alpha}}\right]\right\}^{1/(\beta-\alpha)}. \qquad (10.5)$$

As one can note, the rule of the power of a $q$-exponential was used in eq. (10.3): $\left(e_q^z\right)^r = e_{1-(1-q)/r}^{rz}$. The other solution is obtained in a similar way but starting with $bx^\beta(1+(a/b)x^{\alpha-\beta})= -c$. For example, the solution of $x^\pi + x^e - 1 = 0$ is

$$x = \left\{\frac{\pi}{e-\pi}W_{1-\frac{\pi}{e-\pi}}\left[\frac{e-\pi}{\pi}\right]\right\}^{1/(e-\pi)} \approx 0.7890 \qquad (11)$$

while the solution of $x^{2.01}-5x+6 = 0$ is

$$x = \left\{\frac{2.01}{5.05}W_{1+\frac{2.01}{1.01}}\left[\frac{5.05}{2.01}(-6)^{-\frac{1.01}{2.01}}\right]\right\}^{-1/(1.01)} \approx 2.0302. \qquad (12)$$

As expected, the solution of eq. (12) is close to 2 (solution of $x^2 - 5x + 6$).

The Fermat equation is very famous because of the last Fermat's theorem, solved by Andrew Wiles in 1994. Although Fermat's theorem has been proved to be true, no analytical solution for the Fermat equation has been provided up to now. The solutions of $A^x + B^x = C^x$ are given by

$$x_1 = \frac{1}{\ln(B/A)} \ln\left[\frac{\ln(A/C)}{\ln(B/a)} W_{1-\frac{\ln(A/C)}{\ln(B/A)}}\left(\left(\frac{\ln(A/C)}{\ln(B/A)}\right)^{-1}\right)\right] \qquad (13)$$

$$x_2 = \frac{1}{\ln(A/B)} \ln\left[\frac{\ln(B/C)}{\ln(A/B)} W_{1-\frac{\ln(B/C)}{\ln(A/B)}}\left(\left(\frac{\ln(B/C)}{\ln(A/B)}\right)^{-1}\right)\right] \qquad (14)$$

*Proof.*

$$A^x + B^x = C^x \Rightarrow \left(\frac{A}{C}\right)^x + \left(\frac{B}{C}\right)^x = 1 \Rightarrow e^{\log\left(\frac{A}{C}\right)^x} + e^{\log\left(\frac{B}{C}\right)^x} = e^{x\log\left(\frac{A}{C}\right)} + e^{x\log\left(\frac{B}{C}\right)} = (e^x)^{\log\left(\frac{A}{C}\right)} + (e^x)^{\log\left(\frac{B}{C}\right)} = 1 \Rightarrow (15.1)$$

$$y^{\log\left(\frac{A}{C}\right)} + y^{\log\left(\frac{B}{C}\right)} - 1 = 0. \qquad (15.2)$$

Since eq. (15.2) is written in the form given by eq. (1), it can be solved by eqs. (8) and (9) using $a = 1$, $\alpha = \log(A/C)$, $b = 1$, $\beta = \log(B/C)$ and $c = -1$. Once the value of $y$ in (15.2) was obtained, the values of $x$ in Fermat equation is obtained by taking the natural logarithm of $y$ what gives precisely eqs. (13) and (14). For example, the solution of $4^x + 3^x = 5^x$ is

$$x_1 = \frac{1}{\ln(3/4)} \ln\left[\frac{\ln(4/5)}{\ln(3/a)} W_{1-\frac{\ln(4/5)}{\ln(3/4)}}\left(\frac{\ln(3/4)}{\ln(4/5)}\right)\right] = 2 \qquad (16)$$

while the solution of $22^x + 4^x = 121^x$ is

$$x_2 = \frac{1}{\ln(22/4)} \ln\left[\frac{\ln(4/121)}{\ln(22/4)} W_{1-\frac{\ln(4/121)}{\ln(22/4)}}\left(\left(\frac{\ln(4/121)}{\ln(22/4)}\right)^{-1}\right)\right] \approx 0.2823. \qquad (17)$$

At last, the Fibonacci equation is considered. It is well known that the *n*-th number of the Fibonacci sequence is given by

$$F_n = \frac{1}{\sqrt{5}}\left[\left(\frac{1+\sqrt{5}}{2}\right)^n - \left(\frac{1-\sqrt{5}}{2}\right)^n\right]. \tag{18}$$

Equation (18) tell us that infinite integer numbers (the Fibonacci sequence) can be generated by using the golden rate $\varphi = (1+\sqrt{5})/2$ when the exponent $n$ is also an integer number. However, one can note that using in eq. (18) a variable that can assume real values and using $\bar{\varphi} = (\sqrt{5}-1)/2$ in order to avoid complex values, any integer number can be generated. In fact, the solutions of

$$\varphi^x - \bar{\varphi}^x = y\sqrt{5} \tag{19}$$

are

$$x_1 = \frac{1}{\ln(\bar{\varphi}/\varphi)}\ln\left[\frac{\ln(1/\varphi)}{\ln(\bar{\varphi}/\varphi)}W_{1+\frac{\ln(1/\varphi)}{\ln(\bar{\varphi}/\varphi)}}\left(-\frac{\ln(\bar{\varphi}/\varphi)}{\ln(\varphi)}(y\sqrt{5})^{\frac{\ln(\bar{\varphi}/\varphi)}{\ln(\varphi)}}\right)\right] \quad (y>0) \tag{20}$$

$$x_2 = \frac{1}{\ln(\varphi/\bar{\varphi})}\ln\left[\frac{\ln(1/\bar{\varphi})}{\ln(\varphi/\bar{\varphi})}W_{1+\frac{\ln(1/\bar{\varphi})}{\ln(\varphi/\bar{\varphi})}}\left(-\frac{\ln(\varphi/\bar{\varphi})}{\ln(\bar{\varphi})}(-y\sqrt{5})^{\frac{\ln(\varphi/\bar{\varphi})}{\ln(\bar{\varphi})}}\right)\right] \quad (y<0), \tag{21}$$

One can note that $y = 0$ has no finite solution since $W_q(0) = 0$ for any value of $q$. Figure 3 shows the plot $x$ versus $y$ when the variable $y$ assumes the integer values in the interval [-500,500], while Fig. 4 shows the same in the interval [-20,20].

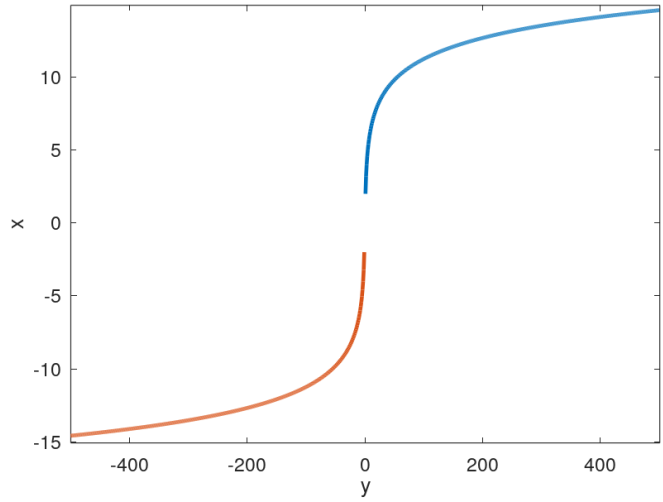

Figure 3 – $x$ versus $y$, eq. (19) – $y \in [-500,500]$ (integer numbers).

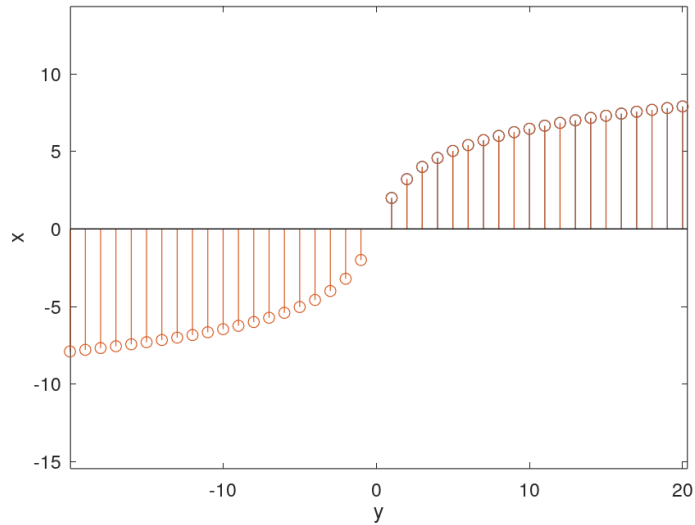

Figure 4 – $x$ versus $y$, eq. (19) – $y \in [-20,20]$ (integer numbers).

One can easily check that eq. (19) can be written in the form shown in eq. (1) ($a = 1$, $\alpha = \varphi$, $b = -1$, $\beta = \bar{\varphi}$ and $c = -y\sqrt{5}$) and, hence, eqs. (8) and (9) can be used to solve eq. (19) ($x = \log(z)$):

$$\varphi^x - \bar{\varphi}^x = y\sqrt{5} \Rightarrow \left(e^x\right)^{\log(\varphi)} - \left(e^x\right)^{\log(\bar{\varphi})} = z^{\log(\varphi)} - z^{\log(\bar{\varphi})} = y\sqrt{5}. \tag{22}$$

## Conclusions

The Lambert-Tsallis $W_q$ function was proposed in 2018 and since them it has been used to solve problems in different areas of physics and engineering like quantum optics, semiconductor physics, plasma physics, quantum information, astronomy, information theory and image processing among others. Here, it was shown the $W_q$ function can also find very interesting applications in pure mathematics. This was shown by providing, for the first time, the general analytical solutions of Fermat and Fibonacci's equations.

## Acknowledgements

This study was financed in part by the Coordenação de Aperfeiçoamento de Pessoal de Nível Superior - Brasil (CAPES) - Finance Code 001, and CNPq via Grant no. 309374/2021-9.